# On non-hamiltonian circulant digraphs of outdegree three

**Stephen C. Locke**
DEPARTMENT OF MATHEMATICAL SCIENCES,
FLORIDA ATLANTIC UNIVERSITY, BOCA RATON, FL 33431

**Dave Witte**
DEPARTMENT OF MATHEMATICS,
OKLAHOMA STATE UNIVERSITY, STILLWATER, OK 74078

**ABSTRACT**

We construct infinitely many connected, circulant digraphs of outdegree three that have no hamiltonian circuit. All of our examples have an even number of vertices, and our examples are of two types: either every vertex in the digraph is adjacent to two diametrically opposite vertices, or every vertex is adjacent to the vertex diametrically opposite to itself.

## 1. INTRODUCTION

It is well known (and not difficult to prove) that every connected, circulant graph has a hamiltonian cycle (except the trivial counterexamples on one or two vertices). (See [2] for much stronger results.) The situation is different in the directed case: some connected, circulant digraphs are not hamiltonian. In general, no good characterization of the hamiltonian circulant digraphs is known. For those of outdegree two, however, R. A. Rankin found a simple arithmetic criterion that determines which are hamiltonian. To state this result, we introduce a bit of notation. (In this paper, circulant digraphs are represented as Cayley digraphs on cyclic groups.)

**Definition 1.1.** For any natural number $n$, we use $\mathbb{Z}_n$ to denote the additive cyclic group of integers modulo $n$. For any set $A$ of integers, let $\mathrm{Cay}(\mathbb{Z}_n; A)$ be the digraph whose vertex set is $\mathbb{Z}_n$, and in which there is an arc from $u$ to $u + a$ (mod $n$), for every $u \in \mathbb{Z}_n$ and every $a \in A$. A digraph is *circulant* if it is (isomorphic to) $\mathrm{Cay}(\mathbb{Z}_n; A)$, for some choice of $n$ and $A$.

Note that $\mathrm{Cay}(\mathbb{Z}_n; A)$ is regular, and its outdegree is equal to the cardinality of the generating set $A$. It is easy to see that $\mathrm{Cay}(\mathbb{Z}_n; A)$ is connected if and only if $\gcd(a_1, a_2, \ldots, a_m, n) = 1$, where $A = \{a_1, a_2, \ldots, a_m\}$.

**Theorem 1.2 (Rankin [5, Thm. 4]).** A connected, circulant digraph $\mathrm{Cay}(\mathbb{Z}_n; a, b)$ of outdegree two has a hamiltonian circuit if and only if there are nonnegative integers $s$ and $t$, such that $s + t = \gcd(sa + tb, n) = \gcd(a - b, n)$.

In contrast, little is known about the hamiltonicity of circulant digraphs of outdegree three (or more). The following theorem provides an interesting class of examples that are hamiltonian.

**Theorem 1.3 (Curran-Witte [4, Thm. 9.1]).** Suppose $\mathrm{Cay}(\mathbb{Z}_n; A)$ is connected, and has outdegree at least three. If $\gcd(a, n) \gcd(b_1, b_2, \ldots, b_m) \geq n$, whenever $a, b_1, b_2, \ldots, b_m \in A$ and $a \notin \{b_1, b_2, \ldots, b_m\}$, then $\mathrm{Cay}(\mathbb{Z}_n; A)$ has a hamiltonian circuit.

| | | | |
|---|---|---|---|
| $\mathrm{Cay}(\mathbb{Z}_{12}; 2, 3, 8)$ | $\mathrm{Cay}(\mathbb{Z}_{28}; 2, 7, 16)$ | $\mathrm{Cay}(\mathbb{Z}_{36}; 2, 9, 20)$ | $\mathrm{Cay}(\mathbb{Z}_{42}; 2, 12, 33)$ |
| $\mathrm{Cay}(\mathbb{Z}_{12}; 3, 4, 6)$ | $\mathrm{Cay}(\mathbb{Z}_{30}; 2, 3, 18)$ | $\mathrm{Cay}(\mathbb{Z}_{36}; 2, 15, 20)$ | $\mathrm{Cay}(\mathbb{Z}_{42}; 2, 15, 36)$ |
| $\mathrm{Cay}(\mathbb{Z}_{18}; 2, 3, 12)$ | $\mathrm{Cay}(\mathbb{Z}_{30}; 2, 6, 21)$ | $\mathrm{Cay}(\mathbb{Z}_{36}; 3, 8, 18)$ | $\mathrm{Cay}(\mathbb{Z}_{42}; 2, 18, 39)$ |
| $\mathrm{Cay}(\mathbb{Z}_{18}; 2, 6, 15)$ | $\mathrm{Cay}(\mathbb{Z}_{30}; 2, 9, 24)$ | $\mathrm{Cay}(\mathbb{Z}_{40}; 2, 5, 22)$ | $\mathrm{Cay}(\mathbb{Z}_{42}; 3, 14, 24)$ |
| $\mathrm{Cay}(\mathbb{Z}_{20}; 2, 5, 12)$ | $\mathrm{Cay}(\mathbb{Z}_{30}; 2, 10, 25)$ | $\mathrm{Cay}(\mathbb{Z}_{40}; 4, 5, 24)$ | $\mathrm{Cay}(\mathbb{Z}_{42}; 6, 7, 28)$ |
| $\mathrm{Cay}(\mathbb{Z}_{24}; 2, 3, 14)$ | $\mathrm{Cay}(\mathbb{Z}_{30}; 3, 10, 18)$ | $\mathrm{Cay}(\mathbb{Z}_{42}; 2, 3, 24)$ | $\mathrm{Cay}(\mathbb{Z}_{44}; 2, 11, 24)$ |
| $\mathrm{Cay}(\mathbb{Z}_{24}; 2, 9, 12)$ | $\mathrm{Cay}(\mathbb{Z}_{30}; 5, 6, 20)$ | $\mathrm{Cay}(\mathbb{Z}_{42}; 2, 6, 27)$ | |
| $\mathrm{Cay}(\mathbb{Z}_{24}; 3, 4, 16)$ | $\mathrm{Cay}(\mathbb{Z}_{36}; 2, 3, 20)$ | $\mathrm{Cay}(\mathbb{Z}_{42}; 2, 7, 28)$ | |

FIGURE 1. Non-hamiltonian, connected, circulant digraphs of outdegree 3 with less than 48 vertices.

One non-hamiltonian example, $\mathrm{Cay}(\mathbb{Z}_{12}; 3, 4, 6)$, was found by D. Witte [6, p. 301]. In this paper, we construct infinitely many non-hamiltonian, connected, circulant digraphs of outdegree three (without loops or multiple arcs). (Figure 1 lists examples with less than 48 vertices. For brevity, the table does not list $\mathrm{Cay}(\mathbb{Z}_n; xa, xb, xc)$ if it includes $\mathrm{Cay}(\mathbb{Z}_n; a, b, c)$, and $\gcd(x, n) = 1$.) In all of our examples, $n$ is even, and the examples come in two families: either the generating set $A$ contains the element $n/2$ of order two in $\mathbb{Z}_n$ (see 3.1'), or two of the elements of $A$ differ by $n/2$ (see 4.6').

**Theorem 3.1'.** For $k \geq 1$, the circulant digraph $\mathrm{Cay}(\mathbb{Z}_{12k}; 6k, 6k+2, 6k+3)$ has no hamiltonian circuit.

If $\gcd(x, n) = 1$, then $\mathrm{Cay}(\mathbb{Z}_n; xa, xb, xc)$ is isomorphic to $\mathrm{Cay}(\mathbb{Z}_n; a, b, c)$, so this theorem can be restated in the following more general form.

**Corollary 1.4.** If $\gcd(a - b, 12k) = 1$, and either $2a - 3b \equiv 6k \pmod{12k}$ or $3a - 2b \equiv 6k \pmod{12k}$, then $\mathrm{Cay}(\mathbb{Z}_{12k}; 6k, a, b)$ has no hamiltonian circuit.

**Theorem 4.6'.** The circulant digraph $\mathrm{Cay}(\mathbb{Z}_{2k}; a, b, b+k)$ has no hamiltonian circuit if and only if $\gcd(a, b, k) \neq 1$, or

- $\gcd(a - b, k) = 1$; and
- $\gcd(a, 2k) \neq 1$; and
- $\gcd(b, k) \neq 1$; and
- either $a$ or $k$ is odd; and
- $a$ is even, or both of $b$ and $k$ are even.

It is natural to ask whether there are any other non-hamiltonian examples. In this vein, an exhaustive computer search reported that every non-hamiltonian, connected, circulant digraph of outdegree three with no more than 95 vertices is described by either Corollary 1.4 or Theorem 4.6'. (If this computer calculation is correct, then Corollary 5.2 implies that if there exists a connected, non-hamiltonian, circulant digraph with outdegree four (or more), then it must have more than 95 vertices.) Perhaps the first question to ask is whether the converse of Corollary 1.4 is true: if $\mathrm{Cay}(\mathbb{Z}_{2n}; n, a, b)$ has no hamiltonian circuit, must it be the case that $n$ is divisible by 6, $\gcd(a - b, 2n) = 1$, and either $2a - 3b$ or $3a - 2b$ is $\equiv n \pmod{2n}$? More fundamental, but also, presumably, more difficult is to determine whether there are any examples with an odd number of vertices, or of outdegree $\geq 4$.

Our results do not provide any counterexamples to the following conjecture.

**Conjecture 1.5 (Curran-Witte [4, p. 74]).** Suppose $\mathrm{Cay}(\mathbb{Z}_n; A)$ is connected, and has outdegree at least three. If, for every proper subset $A'$ of $A$, the subdigraph $\mathrm{Cay}(\mathbb{Z}_n; A')$ is not connected, then $\mathrm{Cay}(\mathbb{Z}_n; A)$ has a hamiltonian circuit. ∎

As mentioned above, circulant digraphs are Cayley digraphs on cyclic groups. Thus, this paper is related to the literature on hamiltonian circuits in Cayley digraphs [1], [3], [6]. Indeed, Rankin's Theorem (1.2) was proved for 2-generated Cayley digraphs on any abelian group, not just on cyclic groups (and even some Cayley digraphs on nonabelian groups). Similarly, Theorem 1.3 and Conjecture 1.5 are only special cases of statements for all abelian groups.

A basic lemma and some definitions are presented in Section 2. The proofs of Theorems 3.1 and 4.6 are given in Sections 3 and 4, respectively. A small result on the hamiltonicity of circulants of outdegree four or more appears in Section 5.

## 2. A PARITY LEMMA

**Definition 2.1.** Given a digraph $G$, let $\mathcal{C} = \mathcal{C}(G)$ be the set of all spanning subdigraphs of $G$ with indegree 1 and outdegree 1 at each vertex. (Thus, each component of a digraph in $\mathcal{C}$ is a circuit.)

**Lemma 2.2.** Given a digraph $G$, suppose $H$ and $H'$ belong to $\mathcal{C}$. Let $u_1$, $u_2$, and $u_3$ be three vertices of $H$, and let $v_i$ be the vertex that follows $u_i$ in $H$. Assume $H'$ has the same arcs as $H$, except:

- instead of the arcs from $u_1$ to $v_1$, from $u_2$ to $v_2$, and from $u_3$ to $v_3$,
- there are arcs from $u_1$ to $v_2$, from $u_2$ to $v_3$, and from $u_3$ to $v_1$.

Then the number of components of $H$ has the same parity as the number of components of $H'$.

*Proof.* Let $\sigma$ be the permutation of $\{1, 2, 3\}$ defined by: $u_{\sigma(i)}$ is the vertex that is encountered when $H$ first reenters $\{u_1, u_2, u_3\}$ after $u_i$. Thus, if $\sigma$ is the identity permutation, then $u_1, u_2, u_3$ lie on three different components of $H$. On the other hand, if $\sigma$ is a 2-cycle, then two of $u_1, u_2, u_3$ are on the same component, but the third is on a different component. Similarly, if $\sigma$ is a 3-cycle, then all three of these vertices are on the same component. Thus, the parity of the number of components of $H$ that intersect $\{u_1, u_2, u_3\}$ is precisely the opposite of the parity of the permutation $\sigma$.

There is a similar permutation $\sigma'$ for $H'$. From the definition of $H'$, we see that $\sigma'$ is simply the product of $\sigma$ with the 3-cycle $(1, 2, 3)$, so $\sigma'$ has the same parity as $\sigma$, because 3-cycles are even permutations. Thus, the parity of the number of components of $H$ that intersect $\{u_1, u_2, u_3\}$ is the same as the parity of the number of components of $H'$ that intersect $\{u_1, u_2, u_3\}$. Because the components that do not intersect $\{u_1, u_2, u_3\}$ are exactly the same in $H$ as in $H'$, this implies that the number of components in $H$ has the same parity as the number of components in $H'$. ∎

**Definition 2.3.** Let $G = \operatorname{Cay}(\mathbb{Z}_n; A)$, and suppose $H \in \mathcal{C}$. For any $u \in \mathbb{Z}_n$ and $a \in A$, we say that $u$ *travels by* $a$ in $H$ if the arc from $u$ to $u + a$ is in $H$.

## 3. A GENERATOR OF ORDER TWO

**Theorem 3.1.** If $a$ is divisible by 6, then $\operatorname{Cay}(\mathbb{Z}_{2a}; a, a+2, a+3)$ has no hamiltonian circuit.

*Proof.* Suppose there is a hamiltonian circuit $H_0$; let $r$ be the number of vertices that travel by $a$, let $s$ be the number of vertices that travel by $a + 2$, and let $t$ be the number of vertices that travel by $a + 3$. Since $a$ and $a + 2$ are both even, we have $\gcd(a, a+2, 2a) \neq 1$, so $t \neq 0$. Also, since $a$ is divisible by 3, we have $\gcd(a+3, 2a) \neq 1$, so $t \neq 2a$. Therefore, $0 < t < 2a$.

We must have $r + s + t = 2a$, and $ra + s(a+2) + t(a+3)$ must be divisible by $2a$. Therefore, we have

$$t = \bigl(ra + s(a+2) + t(a+3)\bigr) - (a+2)(r+s+t) + 2r \equiv 2r \pmod{2a},$$

and

$$s = (a+3)(r+s+t) - \bigl(ra + s(a+2) + t(a+3)\bigr) - 3r \equiv -3r \pmod{2a}.$$

Now $r \leq a$, because the hamiltonian circuit can never have two consecutive $a$-arcs. Therefore, because $0 < t < 2a$, the congruence $t \equiv 2r$ implies that

$$t = 2r.$$

Therefore, we have $2s + 3t = t + 2s + 2t = 2(r + s + t) = 2(2a) = 4a$.

For each $i \in \mathbb{Z}_{2a}$, let

$$B_i = \{i, i+1, i+2, a+i, a+i+1, a+i+2\}.$$

We claim that

for each $i$, the subdigraph of $H_0$ induced by $B_i$ has exactly two arcs.

Consider the walk $\overline{W}$ in $\mathrm{Cay}(\mathbb{Z}_a; 2, 3)$ that results from reducing $H$ modulo $a$, and removing the loops. This walk may be lifted to a path $W$ in $\mathrm{Cay}(\mathbb{Z}; 2, 3)$ that begins at 0 and ends at $4a$. Thus, for each $j$, with $0 \leq j < 4a$, there is exactly one arc $u_j \to v_j$ of $W$ with $u_j \leq j$ and $v_j > j$. Because $j < v_j \in \{u_j + 2, u_j + 3\}$, we have $u_j \geq j - 2$, so the arc $u_j \to v_j$ starts in the set $\{j-2, j-1, j\}$ and ends outside this set. The corresponding arc $\overline{u}_j \to \overline{v}_j$ of $H$ starts in $B_{j-2}$ and ends outside $B_{j-2}$. Because $B_{j-2} = B_i$ iff $j - 2 \equiv i \pmod{a}$, we conclude that the hamiltonian circuit $H_0$ has exactly 4 arcs that start in $B_i$ and end outside $B_i$. The claim follows.

Let $\mathcal{D}$ be the collection of all spanning subdigraphs $H$ of $\mathrm{Cay}(\mathbb{Z}_{2a}; a, a+2, a+3)$, such that

(1) every vertex of $H$ has indegree 1 and outdegree 1 (that is, $H \in \mathcal{C}$);
(2) $H$ has an odd number of components;
(3) we have $t = 2r$, where $t = t_H$ is the number of vertices that travel by $a+3$ in $H$, and $r = r_H$ is the number that travel by $a$; and
(4) for each $i$, the subdigraph of $H$ induced by $B_i$ has exactly two arcs.

We know $\mathcal{D}$ is nonempty, because the hamiltonian circuit $H_0$ belongs to $\mathcal{D}$.

Let $H$ be a digraph in $\mathcal{D}$, such that $r$ is minimal.

We claim that some vertex travels by $a$ in $H$. For, otherwise, we have $r = r_H = 0$, which implies $t = 2r = 0$, so every vertex of $H$ must travel by $a+2$. Therefore, the number of components of $H$ is precisely $\gcd(a+2, 2a)$. Because $a$ is even (indeed, it is divisible by 6), this implies that $H$ has an even number of components, which contradicts the definition of $\mathcal{D}$.

*Case* 1. *For some $i$, the two consecutive vertices $i$ and $i+1$ both travel by $a$ in $H$.* By vertex-transitivity, there is no harm in assuming $i = a+1$. Since the two arcs $(a+1) \to 1$ and $(a+2) \to 2$ must be the only arcs within the blocks $B_1$ and $B_2$, we see that 0, $a$, and 1 must all travel by $a+3$. For the same reason, the vertex 2 cannot travel by $a$. However, the vertex 2 cannot travel by $a+2$, lest the vertex $a+4$ have indegree two; so the vertex 2 must travel by $a+3$. Then the vertex 3 must also travel by $a+3$, lest either $a+3$ or $a+5$ have indegree two. Continuing this argument, we see that 4, 5, 6,... must all travel by $a+3$. So every vertex travels by $a+3$, which contradicts the assumption that $i$ travels by $a$.

*Case* 2. *For every $i$, if the vertex $i$ travels by $a$, then the vertex $i-2$ also travels by $a$.* Some vertex travels by $a$, so, by vertex-transitivity, there is no harm in assuming that 0 travels by $a$. Hence, by repeated application of the hypothesis, we see that the vertex $2j$ travels by $a$, for every $j$. In particular, the vertices 0, 2, $a$, and $a+2$ all travel by $a$. This contradicts the fact that the subdigraph of $H$ induced by the block $B_0$ has only two arcs.

*Case* 3. *The general case.* Some vertex travels by $a$, so, by vertex-transitivity, there is no harm in assuming that 3 travels by $a$. From Case 2, we may assume that 1 does not travel by $a$. However, the vertex 1 also does not travel by $a+2$, lest the vertex $a+3$ have indegree two; thus, the vertex 1 must travel by $a+3$.

Now, from Case 1, we may assume that the vertex 2 does not travel by $a$. However, it also does not travel by $a+2$, lest the vertex $a+4$ have indegree two; hence, the vertex 2 must travel by $a+3$.

Now, we construct another spanning subdigraph $H'$ in which the vertices 1, 2, and 3 all travel by $a+2$: $H'$ has the same arcs as $H$, except:

- instead of the arcs from 1 to $a+4$, from 2 to $a+5$, and from 3 to $a+3$,
- there are arcs from 1 to $a+3$, from 2 to $a+4$, and from 3 to $a+5$.

Note that $t' = t - 2$ and $r' = r - 1$, so $t' = t - 2 = 2r - 2 = 2r'$.

From Lemma 2.2, we know that the number of components of $H$ has the same parity as the number of components of $H'$. That is, the number of components of $H'$ is odd. We conclude that $H' \in \mathcal{D}$. But, because $r' = r - 1$, this contradicts the minimality of $H$. ∎

## 4. GENERATORS WHOSE DIFFERENCE IS THE ELEMENT OF ORDER TWO

**Definition 4.1.** Let $G = \mathrm{Cay}(\mathbb{Z}_{2k}; a, b, b + k)$. Let $\mathcal{E} = \mathcal{E}(G)$ be the set of all spanning subdigraphs of $G$ with indegree 1 and outdegree 1 at each vertex, such that, in each coset of the subgroup $\{0, k\}$, exactly one vertex travels by $a$, and the other by $b$ or $b + k$. (Note that $\mathcal{E}$ is a subset of the class $\mathcal{C}$ introduced in §2.)

**Notation 4.2.** For any subset $A$ of a group $\Gamma$, we use $\langle A \rangle$ to denote the subgroup of $\Gamma$ generated by $A$. For $A \subset \mathbb{Z}_n$, note that $\mathrm{Cay}(\Gamma; A)$ is connected if and only if $\langle A \rangle = \mathbb{Z}_n$.

**Definition 4.3.** Let $G = \mathrm{Cay}(\mathbb{Z}_{2k}; a, b, b + k)$, and assume $G$ is connected. We construct an element $H_0$ of $\mathcal{E}$. Let $d = 2k/\gcd(a, 2k)$ be the order of the element $a$ in the cyclic group $\mathbb{Z}_{2k}$; the construction of our example depends on the parity of $d$.

*Case* 1. *$d$ is odd.* In this case, $k \notin \langle a \rangle$. Every vertex $v$ in $\mathbb{Z}_{2k}$ can be uniquely written in the form $x_v a + y_v b + z_v k$ with $0 \leq x_v < d$, $0 \leq y_v < k/d$, and $0 \leq z_v < 2$. Let $H_0$ be the spanning subdigraph in which a vertex $v \in \mathbb{Z}_{2k}$

- travels by $a$ if $z_v = 0$;
- travels by $b$ if $z_v = 1$ and $z_{v+b} = 1$; and
- travels by $b + k$ otherwise.

(By construction, the vertices $v$ that satisfy $z_v = 0$ are both entered and exited via an $a$-arc in $H_0$; the other vertices are neither entered nor exited via an $a$-arc.)

*Case* 2. *$d$ is even.* In this case, $k \in \langle a \rangle$, so every vertex $v$ in $\mathbb{Z}_{2k}$ can be uniquely written in the form $x_v a + y_v b$ with $0 \leq x_v < d$ and $0 \leq y_v < 2k/d$. Let $H_0$ be the spanning subdigraph in which a vertex $v \in \mathbb{Z}_{2k}$

- travels by $a$ if $x_v < d/2$;
- travels by $b + k$ if $x_v \geq d/2$ and $1 \leq x_{v+b} \leq d/2$; and
- travels by $b$ otherwise.

(By construction, the vertices $v$ that satisfy $1 \leq x_v \leq d/2$ are precisely those that are entered via an $a$-arc in $H_0$.)

**Lemma 4.4.** Let $G = \mathrm{Cay}(\mathbb{Z}_{2k}; a, b, b + k)$, assume $G$ is connected, and let $H_0$ be the element of $\mathcal{E}$ constructed in Definition 4.3. Then $H_0$ has an odd number of components if and only if either

- both of $a$ and $k$ are even; or
- $a$ is odd, and either $b$ or $k$ is odd.

*Proof.* Let $d = 2k/\gcd(a, 2k)$ be the order of the element $a$ in the cyclic group $\mathbb{Z}_{2k}$; the proof depends on the parity of $d$.

*Case* 1. *$d$ is odd.* Because $ad$ is a multiple of $2k$, we see, in this case, that $a$ must be even. Thus, we wish to show that the parity of the number of components of $H_0$ is the opposite of the parity of $k$.

For $i \in \{0, 1\}$, let $G_i = \{ v \in \mathbb{Z}_{2k} \mid z_v = i \}$, so each of $G_0$ and $G_1$ has exactly $k$ vertices. From the definition of $H_0$, we see that each component of $H_0$ is contained in either $G_0$ or $G_1$. Each component in $G_0$ is a circuit of length $d$ (all $a$-arcs), so the number of components in $G_0$ is $k/d$. Because $d$ is odd, this has the same parity as $k$, so we wish to show that $G_1$ contains an odd number of components of $H_0$.

The number of components contained in $G_1$ is equal to the order of the quotient group $\mathbb{Z}_{2k}/\langle b, k \rangle$. Because $\langle a, b, k \rangle = \mathbb{Z}_{2k}$, we know that $a$ generates this quotient group. Then, because $a$ has odd order, we conclude that the quotient group also has odd order, as desired.

*Case* 2. *d is even.* Let $xa + yb$ be a vertex that travels by $a$ in $H_0$. Then $v = (d/2)a + yb$ is in the same component (by following a sequence of $a$-arcs). Furthermore, if $y < (2k/d) - 1$, then we see that $x_{v+b} = d/2$, so $v$ travels by $b + k$; this means that $(y+1)b = v + b + k$ is also in the same component. By induction on $y$, this implies that all the $a$-arcs of $H_0$ are in the same component, and this component contains some $(b + k)$-arcs. Thus, the $a$-arcs are essentially irrelevant in counting components of $H_0$: there is a natural one-to-one correspondence between the components of $H_0$ and the components of $\mathrm{Cay}(\mathbb{Z}_k; b)$. Thus, the number of components is equal to the order of the quotient group $\mathbb{Z}_k/\langle b \rangle$. This quotient group has odd order if and only if either $b$ or $k$ is odd. Therefore, $H_0$ has an odd number of components if and only if either $b$ or $k$ is odd.

Thus, we have the desired conclusion if $a$ is odd, so we may now assume $a$ is even. Since $2k/\gcd(2k, a) = d$ is even, this implies that $k$ is also even. So we wish to show that $H_0$ has an odd number of components. Because $\mathrm{Cay}(\mathbb{Z}_{2k}; a, b, k)$ is connected, it cannot be the case that $a$, $b$, and $k$ are all even, so we conclude that $b$ is odd. From the conclusion of the preceding paragraph, we see that $H_0$ has an odd number of components, as desired. ∎

The following result is a consequence of the proof of Lemma 2.2.

**Lemma 4.5.** Let $G = \mathrm{Cay}(\mathbb{Z}_{2k}; a, b, b+k)$, assume $H \in \mathcal{E}$, and suppose $u$ is a vertex of $H$ that travels by $a$, such that $u$, $u+k$, and $u+a+k$ are on three different components of $H$. Then there is an element $H'$ of $\mathcal{E}$, with exactly the same arcs as $H$, except the arcs leaving $u$ and $u+k$, and the arc entering $u+a+k$, such that $u$, $u+k$, and $u+a+k$ are all on the same component of $H'$.

**Theorem 4.6.** The circulant digraph $\mathrm{Cay}(\mathbb{Z}_{2k}; a, b, b+k)$ has a hamiltonian circuit if and only if $\gcd(a, b, k) = 1$, and either

- $\gcd(a - b, k) \neq 1$; or
- $\gcd(a, 2k) = 1$; or
- $\gcd(b, k) = 1$; or
- both of $a$ and $k$ are even; or
- $a$ is odd, and either $b$ or $k$ is odd.

*Proof.* ($\Rightarrow$) Because hamiltonian digraphs are connected, we know that $\gcd(a, b, k) = 1$. We may assume $\gcd(a - b, k) = 1$, $\gcd(a, 2k) \neq 1$, and $\gcd(b, k) \neq 1$.

Choose a hamiltonian circuit; let $r$ be the number of vertices that travel by $a$, and let $s$ be the number of vertices that travel by $b$ or $b + k$. We must have $r + s = 2k$, and $ra + sb$ must be divisible by $k$. Therefore, we conclude that $r(a - b)$ is divisible by $k$. Since $\gcd(a - b, k) = 1$, this implies $r$ is divisible by $k$. Because $0 \leq r \leq 2k$, this implies $r \in \{0, k, 2k\}$. Because $\gcd(a, 2k) \neq 1$, we know $\langle a \rangle \neq \mathbb{Z}_{2k}$, so we cannot have $r = 2k$; because $\gcd(b, k) \neq 1$, we know $\langle b, k \rangle \neq \mathbb{Z}_{2k}$, so we cannot have $r = 0$. Therefore, we must have $r = k$. So exactly half of the vertices travel by $a$, and the other half travel by $b$ or $b + k$.

Let us show that every hamiltonian circuit belongs to $\mathcal{E}$. That is, in each coset of the subgroup $\{0, k\}$, exactly one vertex travels by $a$, and the other by $b$ or $b + k$. If not, then, from the conclusion of the preceding paragraph, there must be some coset $i + \{0, k\}$ in which both vertices travel by $a$. Therefore, both vertices of $i + a + \{0, k\}$ are entered via $a$, which means that neither of the vertices in $i + a - b + \{0, k\}$ travels by $b$ or $b + k$, so they both must travel by $a$. Repeating the argument, we see that both of the vertices in $i + j(a - b) + \{0, k\}$ travel by $a$, for all $j$. Because $\gcd(a - b, k) = 1$, every vertex in the digraph is of the form $i + j(a - b)$ or $i + j(a - b) + k$, so we see that every vertex travels by $a$. This contradicts the conclusion of the preceding paragraph.

Recall the digraph $H_0$ of Definition 4.3. It suffices to show, for every $H \in \mathcal{E}$, that the number of components of $H$ has the same parity as the number of components of $H_0$. For then, because the preceding paragraph implies that $\mathcal{E}$ contains a hamiltonian circuit, we conclude that $H_0$ has an odd number of components. Then Lemma 4.4 provides the desired conclusion.

Let $u_1$ be some vertex that travels by $a$ in $H$, and let $v_1 = u_1 + a$. Let $u_2 = u_1 + k$, and let $v_2 \in u_2 + \{b, b+k\}$ be the vertex that follows $u_2$ in $H$. Finally, let $v_3 = v_1 + k$, and let $u_3 \in v_3 - \{b, b+k\}$ be the vertex that *precedes* $v_3$ in $H$. We construct an element $H'$ of $\mathcal{E}$ in which it is $u_2$ that travels by $a$, instead of $u_1$: $H'$ has the same arcs as $H$, except:

- instead of the arcs from $u_1$ to $v_1$, from $u_2$ to $v_2$, and from $u_3$ to $v_3$,
- there are arcs from $u_1$ to $v_2$, from $u_2$ to $v_3$, and from $u_3$ to $v_1$.

Lemma 2.2 implies that the number of components of $H$ has the same parity as the number of components of $H'$.

Because $H$ and $H_0$ both have the property that, in each coset of $\{0, k\}$, exactly one vertex travels by $a$, and the other by $b$ or $b + k$, we may transform $H$ into $H_0$, by performing a sequence of transformations of the form $H \mapsto H'$. Thus, we may transform $H$ into $H_0$, without changing the parity of the number of components, as desired.

($\Leftarrow$) Because $\gcd(a, b, k) = 1$, we know that $\langle a, b, k \rangle = \mathbb{Z}_{2k}$.

*Case* 1. *We have* $\gcd(a, 2k) = 1$. In this case, we have $\langle a \rangle = \mathbb{Z}_{2k}$, so there is an obvious hamiltonian circuit in the Cayley digraph (all $a$-arcs).

*Case* 2. *We have* $\gcd(b, k) = 1$ In this case, either $\langle b \rangle = \mathbb{Z}_{2k}$ or $\langle b+k \rangle = \mathbb{Z}_{2k}$, so there is again an obvious hamiltonian circuit.

*Case* 3. *We have* $\gcd(a - b, k) \neq 1$. In this case, we have $\langle a - b, k \rangle \neq \mathbb{Z}_{2k}$. There are many digraphs in $\mathcal{C}$ in which

- every vertex not in $\langle a - b, k \rangle$ travels by either $b$ or $b + k$; and
- for each vertex $v \in \langle a - b, k \rangle$, one of $v$ and $v + k$ travels by $a$, and the other travels by either $b$ or $b + k$.

Among all such digraphs, let $H$ be one in which the number of components is minimal.

We claim that $H$ is a hamiltonian circuit. If not, then $H$ has more than one component. Because $\langle a, b, k \rangle = \mathbb{Z}_{2k}$, we know that $b$ generates the quotient group $\mathbb{Z}_{2k}/\langle a - b, k \rangle$, so every component of $H$ intersects $\langle a - b, k \rangle$, and hence either

- there is some vertex $u$ in $\langle a - b, k \rangle$ such that $u$ and $u + k$ are in different components of $H$; or
- for all $v \in \langle a - b, k \rangle$, the vertices $v$ and $v + k$ are in the same component of $H$, but there is some vertex $u$ in $\langle a - b, k \rangle$ such that $u$ and $u + (a - b)$ are in different components of $H$.

In either case, let $u_1$ be the one of $u$ and $u + k$ that travels by $a$.

Let $v_1 = u_1 + a$. Let $u_2 = u_1 + k$, and let $v_2 \in u_2 + \{b, b + k\}$ be the vertex that follows $u_2$ in $H$. Finally, let $v_3 = v_1 + k$, and let $u_3 \in v_3 - \{b, b + k\}$ be the vertex that *precedes* $v_3$ in $H$. The choice of $u_1$ implies that $u_1$, $u_2$ and $u_3$ do not all belong to the same component of $H$.

Let $w_1$ and $w_2$ be the vertices that *precede* $u_1$ and $u_2$, respectively, on $H$. (So $w_1 = w_2 + k$.)

Let $\sigma$ be the permutation of $\{1, 2, 3\}$ defined in the proof of Lemma 2.2. If $\sigma$ is an even permutation, let $H_1 = H$; if $\sigma$ is an odd permutation, let $H_1$ be the element of $\mathcal{C}$ that has the same arcs as $H$, except:

- instead of the arcs from $w_1$ to $u_1$, and from $w_2$ to $u_2$,
- there are arcs from $w_1$ to $u_2$, and from $w_2$ to $u_1$.

In either case, the permutation $\sigma_1$ for $H_1$ is even. Thus, $\sigma_1$ is either trivial or a 3-cycle. If it is a 3-cycle, then $u_1$, $u_2$ and $u_3$ are all contained in a single component of $H_1$, so $H_1$ has less components than $H$, which contradicts the minimality of $H$. Thus, $\sigma_1$ is trivial.

Let $H'$ be the element of $\mathcal{C}$ that has the same arcs as $H_1$, except:

- instead of the arcs from $u_1$ to $v_1$, from $u_2$ to $v_2$, and from $u_3$ to $v_3$,
- there are arcs from $u_1$ to $v_2$, from $u_2$ to $v_3$, and from $u_3$ to $v_1$.

Because $\sigma_1$ is trivial, we see that the permutation $\sigma'$ for $H'$ is the 3-cycle $(1, 2, 3)$. Hence, $u_1$, $u_2$ and $u_3$ are all contained in a single component of $H'$, so $H'$ has less components than $H$, which contradicts the minimality of $H$.

*Case* 4. *Either both of $a$ and $k$ are even; or $a$ is odd, and either $b$ or $k$ is odd.* In this case, Lemma 4.4 asserts that the digraph $H_0$ of Definition 4.3 has an odd number of components. We construct a hamiltonian circuit by amalgamating all of these components into one component. We start with the component containing 0, and use Lemma 4.5 to add the other components to it two at a time.

Note that the assumption of the present case, together with the fact that $\gcd(a, b, k) = 1$, implies that $\gcd(b, k)$ is odd. Furthermore, we may assume that $\gcd(b, k) \neq 1$, for, otherwise, Case 2 applies. Thus, $\gcd(b, k) \geq 3$.

Let $d = 2k/\gcd(a, 2k)$ be the order of the element $a$ in the cyclic group $\mathbb{Z}_{2k}$; the proof depends on the parity of $d$.

*Subcase* 4.1. *$d$ is odd.* Note that two vertices $u$ and $v$ are in the same component of $H_0$ if and only if either

- $z_u = z_v = 0$ and $y_u = y_v$; or
- $z_u = z_v = 1$ and $x_u \equiv x_v \pmod{\gcd(b, k)}$.

Lemma 4.5 implies there is an element $H_0'$ of $\mathcal{E}$, such that $0$, $k$, and $a+k$ are all in the same component of $H_0'$. (The other components of $H_0'$ are components of $H_0$.)

Then Lemma 4.5 implies there is an element $H_1 = (H_0')'$ of $\mathcal{E}$, such that $a+b$, $a+b+k$, and $2a+b+k$ are all in the same component of $H_1$. (The other components of $H_1$ are components of $H_0$.)

With this as the base case of an inductive construction, we construct, for $1 \leq i \leq k/(2d)$, an element $H_i$ of $\mathcal{E}$, such that

$$\{\, v \mid z_v = 0 \text{ and } 0 \leq y_v \leq 2i-1 \,\} \cup \{\, v \mid z_v = 1 \text{ and } x_v \equiv 0, 1, \text{ or } 2 \pmod{\gcd(b,k)} \,\}$$

is a component of $H_i$, and all other components of $H_i$ are components of $H_0$. Namely, $H_i$ has exactly the same arcs as $H_{i-1}$, except:

- instead of the arcs

$$\begin{aligned}
(2i-2)b &\to a + (2i-2)b \\
(2i-2)b + k &\to (2i-1)b + k \\
a + (2i-3)b + k &\to a + (2i-2)b + k
\end{aligned}$$

$$\begin{aligned}
(2i-1)b &\to a + (2i-1)b \\
(2i-1)b + k &\to v \\
a + (2i-2)b + k &\to a + (2i-1)b + k
\end{aligned}$$

(where $v = (2i)b + k$ if $i < k/(2d)$, and $v \in \{(2i)b, (2i)b+k\}$ if $i = k/(2d)$),

- there are arcs

$$\begin{aligned}
(2i-2)b &\to (2i-1)b + k \\
(2i-2)b + k &\to a + (2i-2)b + k \\
a + (2i-3)b + k &\to a + (2i-2)b
\end{aligned}$$

$$\begin{aligned}
(2i-1)b &\to v \\
(2i-1)b + k &\to a + (2i-1)b + k \\
a + (2i-2)b + k &\to a + (2i-1)b
\end{aligned}$$

Let $K_1 = H_{k/(2d)}$. With this as the base case of an inductive construction, we construct, for $1 \leq i \leq (\gcd(b,k) - 1)/2$, an element $K_i$ of $\mathcal{E}$, such that

$$\{\, v \mid z_v = 0 \,\} \cup \{\, v \mid z_v = 1 \text{ and } x_v \equiv 0, 1, \ldots, \text{ or } 2i \pmod{\gcd(b,k)} \,\}$$

is a component of $K_i$, and all other components of $K_i$ are components of $H_0$. Namely, Lemma 4.5 implies there is an element $K_i = K_{i-1}'$ of $\mathcal{E}$, such that $(2i-1)a$, $(2i-1)a+k$, and $(2i)a+k$ are all in the same component of $K_i$.

Then, for $i = \bigl(\gcd(b,k)-1\bigr)/2$, we see that a single component of $K_i$ contains every vertex, so $K_i$ is a hamiltonian circuit.

*Subcase* 4.2. *$d$ is even.* Note that one component of $H_0$ is

$$\{\, v \mid x_v < d/2 \,\} \cup \{\, v \mid x_v \equiv 0 \pmod{\gcd(b,k)} \,\}.$$

Two vertices $u$ and $v$ that are not in this component are in the same component of $H_0$ if and only if $x_u \equiv x_v \pmod{\gcd(b,k)}$.

We may assume $2k/d > 1$, for otherwise Case 1 applies. With $H_0$ as the base case of an inductive construction, we construct, for $0 \leq i \leq \bigl(\gcd(b,k) - 1\bigr)/2$, an element $H_i$ of $\mathcal{E}$, such that

$$\{\, v \mid x_v < d/2 \,\} \cup \{\, v \mid x_v \equiv 0, 1, \ldots, \text{ or } 2i \pmod{\gcd(b,k)} \,\}$$

is a component of $H_i$, and all other components of $H_i$ are components of $H_0$. Namely, Lemma 4.5 implies there is an element $H_i = H_{i-1}'$ of $\mathcal{E}$, such that $(2i-1)a$, $(2i-1)a+k$, and $(2i)a+k$ are all in the same component of $H_i$.

Then, for $i = \bigl(\gcd(b,k)-1\bigr)/2$, we see that a single component of $H_i$ contains every vertex, so $H_i$ is a hamiltonian circuit. ∎

## 5. OUTDEGREE AT LEAST FOUR

**Proposition 5.1.** Suppose $\mathrm{Cay}(\mathbb{Z}_n; A)$ has outdegree four or more, and assume there is a proper subset $A'$ of $A$, such that $\mathrm{Cay}(\mathbb{Z}_n; A')$ is connected and has outdegree three. If every non-hamiltonian, connected, circulant digraph that has outdegree three and exactly $n$ vertices is described by either Corollary 1.4 or Theorem 4.6′, then $\mathrm{Cay}(\mathbb{Z}_n; A)$ has a hamiltonian circuit.

*Proof.* Suppose the contrary. Then the spanning subdigraph $\mathrm{Cay}(\mathbb{Z}_n; A')$ also has no hamiltonian circuit. Therefore, by assumption, there are two cases to consider.

*Case* 1. $\mathrm{Cay}(\mathbb{Z}_n; A')$ *is described by Corollary 1.4.* We have $n = 12k$, and there is no harm in assuming that $\mathrm{Cay}(\mathbb{Z}_n; A')$ is described by Theorem 3.1′, so $A' = \{6k, a, b\}$, where $a = 6k + 2$ and $b = 6k + 3$. Let $c$ be an element of $A$ that is not in $A'$. Because $6k \notin \{a, b, c\}$, we know that $\mathrm{Cay}(\mathbb{Z}_n; a, b, c)$ is not described by Corollary 1.4, so it must be described by Theorem 4.6′. Thus, we must have $c \in \{a + 6k, b + 6k\}$. Because both of $a$ and $6k$ are even, we see from Theorem 4.6 that $\mathrm{Cay}(\mathbb{Z}_{12k}; a, b, b + 6k)$ has a hamiltonian circuit. Therefore, it must be the case that $c = a + 6k \equiv 2 \pmod{n}$, so $\{2, 6k, 6k + 2, 6k + 3\} \subset A$. Let $H$ be the spanning subdigraph of $\mathrm{Cay}(\mathbb{Z}_n; A)$ in which every vertex travels by 2, except:

- the vertex 2 travels by $6k$;
- the vertex $6k$ travels by $6k + 2$; and
- the vertices 0 and $6k + 1$ travel by $6k + 3$.

Then $H$ is a hamiltonian circuit.

*Case* 2. $\mathrm{Cay}(\mathbb{Z}_n; A')$ *is described by Theorem 4.6′.* Writing $n = 2k$, we have $A' = \{a, b, b + k\}$; let $c$ be an element of $A$ that is not in $A'$. By interchanging $b$ and $b + k$ if necessary, we may assume $\mathrm{Cay}(\mathbb{Z}_n; a, b)$ is connected. Then we may assume $\mathrm{Cay}(\mathbb{Z}_n; a, b, c)$ is described by Theorem 4.6′, for otherwise Case 1 applies. Therefore, $c \in \{a+k, b+k\}$, so, because $c \notin A'$, we must have $c = a + k$. Any Euler circuit in $\mathrm{Cay}(\mathbb{Z}_k; a, b)$ passes through each vertex exactly twice; any such circuit may be lifted to a hamiltonian circuit in $\mathrm{Cay}(\mathbb{Z}_{2k}; a, a + k, b, b + k)$, which is a spanning subdigraph of $\mathrm{Cay}(\mathbb{Z}_n; A)$. ∎

**Corollary 5.2.** Suppose $\mathrm{Cay}(\mathbb{Z}_n; A)$ is connected, and has outdegree four or more, and assume $n < 420$. If every non-hamiltonian, connected, circulant digraph that has outdegree three and exactly $n$ vertices is described by either Corollary 1.4 or Theorem 4.6′, then $\mathrm{Cay}(\mathbb{Z}_n; A)$ has a hamiltonian circuit.

*Proof.* From the proposition, we may assume there is no 3-element subset $\{a, b, c\}$ of $A$ with $\gcd(a, b, c, n) = 1$. This implies that $n$ has at least four distinct prime factors. Then, since $n < 420 = 2^2 \cdot 3 \cdot 5 \cdot 7$, we know that $n$ is square free. Therefore, because $n < 2310 = 2 \cdot 3 \cdot 5 \cdot 7 \cdot 11$, this implies that $n$ is the product of four distinct primes. Hence, there are four elements $\{a, b, c, d\}$ of $A$ with $\gcd(a, b, c, d, n) = 1$, so we may assume that $A$ has exactly four elements. These conditions imply that the hypotheses of Theorem 1.3 are satisfied, so $\mathrm{Cay}(\mathbb{Z}_n; A)$ has a hamiltonian circuit. ∎

*Remark.* The proof of Corollary 5.2 is much simpler (namely, the first two sentences suffice) if $n < 2 \cdot 3 \cdot 5 \cdot 7 = 210$.

## ACKNOWLEDGMENTS

Witte was partially supported by a grant from the National Science Foundation.